\documentclass[12pt]{article}

\usepackage{amsmath, epsfig, cite}
\usepackage{amssymb}
\usepackage{amsfonts}
\usepackage{latexsym}
\usepackage{float}
\usepackage{color}

\newtheorem{thm}{Theorem}[section]

\newtheorem{lem}[thm]{Lemma}

\numberwithin{equation}{section}

\newcommand{\qed}{{\hfill$\square$}\medskip}

\begin{document}

\begin{center}
{\Large\bf Supercongruences for sums involving Domb numbers}
\end{center}

\vskip 2mm \centerline{Ji-Cai Liu}
\begin{center}
{\footnotesize Department of Mathematics, Wenzhou University, Wenzhou 325035, PR China\\
{\tt jcliu2016@gmail.com } \\[10pt]
}
\end{center}


\vskip 0.7cm \noindent{\bf Abstract.}
We prove some supercongruence and divisibility results on sums involving Domb numbers, which confirm four conjectures of Z.-W. Sun and Z.-H. Sun. For instance, by using a transformation formula due to Chan and Zudilin, we show that for any prime $p\ge 5$,
\begin{align*}
\sum_{k=0}^{p-1}\frac{3k+1}{(-32)^k}{\rm Domb}(k)\equiv (-1)^{\frac{p-1}{2}}p+p^3E_{p-3}
\pmod{p^4},
\end{align*}
which is regarded as a $p$-adic analogue of the following interesting formula for $1/\pi$ due to Rogers:
\begin{align*}
\sum_{k=0}^{\infty}\frac{3k+1}{(-32)^k}{\rm Domb}(k)=\frac{2}{\pi}.
\end{align*}
Here ${\rm Domb}(n)$ and $E_n$ are the famous Domb numbers and Euler numbers.


\vskip 2mm
\noindent{\it MR Subject Classifications}: 11A07, 11Y55, 05A19, 33F10

\section{Introduction}
In 1960, Domb \cite{domb-ap-1960} first introduced the following sequence:
\begin{align*}
{\rm Domb}(n)=\sum_{k=0}^n{n\choose k}^2{2k\choose k}{2n-2k\choose n-k},
\end{align*}
which are known as the famous Domb numbers. This sequence plays an important role in many research fields, including probability theory \cite{bnsw-rj-2011}, special functions \cite{bbbg-2008-jpa}, Ap\'ery-like differential equations \cite{avz-pems-2011}, and combinatorics \cite{rs-ejc-2009}.

The Domb numbers are also connected to some interesting series for $1/\pi$. For instance,
Chan, Chan and Liu \cite{ccl-am-2004} showed that
\begin{align*}
\sum_{k=0}^{\infty}\frac{5k+1}{64^k}{\rm Domb}(k)=\frac{8}{\sqrt{3}\pi}.
\end{align*}
Another typical example is the the following identity due to Rogers \cite{rogers-rj-2009}:
\begin{align}
\sum_{k=0}^{\infty}\frac{3k+1}{(-32)^k}{\rm Domb}(k)=\frac{2}{\pi}.\label{a-rogers}
\end{align}

Let $E_n$ denote the Euler numbers given by
\begin{align*}
\frac{2}{e^x+e^{-x}}=\sum_{n=0}^{\infty}E_n\frac{x^n}{n!}.
\end{align*}
The motivation of this paper is to prove the following interesting $p$-adic analogue of \eqref{a-rogers}, which was originally conjectured by Z.-W. Sun \cite[Conjecture 77 (ii)]{sunzw-2019}.
\begin{thm}\label{t-1}
For any prime $p\ge 5$, we have
\begin{align}
\sum_{k=0}^{p-1}\frac{3k+1}{(-32)^k}{\rm Domb}(k)\equiv (-1)^{\frac{p-1}{2}}p+p^3E_{p-3}
\pmod{p^4}.\label{a-thm1}
\end{align}
\end{thm}

The proof of \eqref{a-thm1} heavily relies on the transformation formula due to Chan and Zudilin \cite[Corollary 3.4]{cz-m-2010}:
\begin{align}
\sum_{k=0}^n{n\choose k}^2{2k\choose k}{2n-2k\choose n-k}=
\sum_{k=0}^n(-1)^k{n+2k\choose 3k}{2k\choose k}^2{3k\choose k}16^{n-k}.\label{a-cz}
\end{align}

The second purpose of this paper is to prove a related supercongruence conjectured by Z.-H. Sun \cite[Conjecture 2.6]{sunzh-jdea-2018} and two divisibility results on sums of Domb numbers conjectured by Z.-W. Sun \cite[Conjecture 77 (i)]{sunzw-2019}.
\begin{thm}\label{t-2}
For any prime $p\ge 5$, we have
\begin{align}
\sum_{k=0}^{p-1}\frac{3k+2}{(-2)^k}{\rm Domb}(k)\equiv 2p(-1)^{\frac{p-1}{2}}+6p^3E_{p-3}
\pmod{p^4}.\label{a-thm2}
\end{align}
\end{thm}
We remark that Z.-W. Sun \cite{sunzw-b-2013} conjectured the supercongruence \eqref{a-thm2} modulo $p^3$.

\begin{thm}\label{t-3}
Let $n$ be a positive integer. Then
\begin{align*}
\frac{1}{n}\sum_{k=0}^{n-1}(2k+1){\rm Domb}(k)8^{n-1-k}
\quad \text{and}\quad\frac{1}{n}\sum_{k=0}^{n-1}(2k+1){\rm Domb}(k)(-8)^{n-1-k}
\end{align*}
are all positive integers.
\end{thm}

The sums of cubes of binomial coefficients:
\begin{align*}
f_n=\sum_{k=0}^n{n\choose k}^3
\end{align*}
are known as Franel numbers \cite{franel-1984}.
The proofs of Theorems \ref{t-2} and \ref{t-3} respectively make use of
the identity due to Z.-H. Sun \cite[Lemma 3.1]{sunzh-itsf-2015}:
\begin{align}
\sum_{k=0}^n{n\choose k}^2{2k\choose k}{2n-2k\choose n-k}=\sum_{k=0}^{\lfloor n/2\rfloor}{2k\choose k}^2{3k\choose k}
{n+k\choose 3k}4^{n-2k},\label{a-sunzh}
\end{align}
and the other identity due to Chan, Tanigawa, Yang and Zudilin \cite[(2.27)]{ctyz-am-2011}:
\begin{align}
\sum_{k=0}^n{n\choose k}^2{2k\choose k}{2n-2k\choose n-k}=(-1)^n\sum_{k=0}^n{n\choose k}{n+k\choose k}(-8)^{n-k}f_k.\label{a-ctyz}
\end{align}

In the past few years, supercongruences for sums of Domb numbers have been widely discussed by many researchers (see, for example, \cite{mw-ijnt-2019,ms-ijnt-2018,sunzh-itsf-2015,sunzh-jdea-2018,sunzw-b-2013,zudilin-jnt-2009}).

The rest of the paper is organized as follows. Section 2 lays down some preparatory results
on combinatorial identities involving harmonic numbers and related congruences.
We prove Theorems \ref{t-1}--\ref{t-3} in Sections 3--5, respectively.

\section{Preliminary results}
Let
\begin{align*}
H_n^{(r)}=\sum_{j=1}^n\frac{1}{j^r}
\end{align*}
denote the $n$th generalized harmonic number of order $r$ with the convention that $H_n=H_n^{(1)}$. The Fermat quotient of an integer $a$ with respect to an odd prime $p$ is given by $q_p(a)=(a^{p-1}-1)/p$.

\begin{lem}
For any non-negative integer $n$, we have
\begin{align}
&\sum_{i=0}^n(-1)^i{n\choose i}{n+i\choose i}\left(H_{2i}-H_i\right)=(-1)^{n+1}\sum_{i=1}^n\frac{(-1)^i}{i},\label{b-1}\\[10pt]
&\sum_{i=0}^n(-1)^i{n\choose i}{n+i\choose i}\left((H_{2i}-H_i)^2-H_{2i}^{(2)}-H_i^{(2)}\right)\notag\\[10pt]
&=2(-1)^n\left(\sum_{i=1}^n\frac{(-1)^i}{i^2}+\sum_{i=1}^n\frac{(-1)^i}{i}H_i\right).
\label{b-2}
\end{align}
\end{lem}
{\noindent\it Proof.}
The identities \eqref{b-1} and \eqref{b-2} are discovered and proved by the symbolic summation package {\tt Sigma} developed by Schneider \cite{schneider-slc-2007}. One can also refer to \cite{liu-jsc-2019,liu-jmaa-2020} for the same approach to finding and proving identities of this type.
\qed

\begin{lem}(See \cite[Lemma 2.4]{sunzw-scm-2011} and \cite[Lemma 2.9]{at-2019}.)
For any prime $p\ge 5$, we have
\begin{align}
&\sum_{i=1}^{(p-1)/2}\frac{(-1)^i}{i^2}\equiv(-1)^{\frac{p-1}{2}}2E_{p-3}\pmod{p},\label{b-3}\\[10pt]
&\sum_{i=1}^{(p-1)/2}\frac{(-1)^i}{i}H_i\equiv \frac{1}{2}q_p(2)^2+(-1)^{\frac{p-1}{2}}E_{p-3}\pmod{p}.\label{b-4}
\end{align}
\end{lem}

\begin{lem}
For any prime $p\ge 5$, we have
\begin{align}
\sum_{i=1}^{(p-1)/2}\frac{(-1)^i}{i}\equiv-q_p(2)+\frac{1}{2}pq_p(2)^2-p(-1)^{\frac{p-1}{2}}E_{p-3}\pmod{p^2}.
\label{b-5}
\end{align}
\end{lem}
{\noindent\it Proof.} We begin with the following congruence \cite[(43)]{lehmer-am-1938}:
\begin{align}
\sum_{i=1}^{\lfloor p/4\rfloor}\frac{1}{p-4i}\equiv \frac{3}{4}q_p(2)-\frac{3}{8}pq_p(2)^2\pmod{p^2}.\label{b-6}
\end{align}
Since for $1\le i\le \lfloor p/4\rfloor$,
\begin{align*}
\frac{1}{p-4i}\equiv -\frac{1}{4i}-\frac{p}{(4i)^2}\pmod{p^2},
\end{align*}
we have
\begin{align}
\sum_{i=1}^{\lfloor p/4\rfloor}\frac{1}{p-4i}\equiv -\frac{1}{4}H_{\lfloor p/4\rfloor}
-\frac{p}{16}H_{\lfloor p/4\rfloor}^{(2)}\pmod{p^2}.\label{b-7}
\end{align}
By \cite[page 359]{lehmer-am-1938}, we have
\begin{align}
H_{\lfloor p/4\rfloor}^{(2)}\equiv (-1)^{\frac{p-1}{2}}4E_{p-3}\pmod{p}. \label{b-8}
\end{align}
Combining \eqref{b-6}--\eqref{b-8}, we arrive at
\begin{align}
H_{\lfloor p/4\rfloor}\equiv -3q_p(2)+\frac{3}{2}pq_p(2)^2-p(-1)^{\frac{p-1}{2}}E_{p-3}
\pmod{p^2}.\label{b-9}
\end{align}

Furthermore, we have
\begin{align}
\sum_{i=1}^{(p-1)/2}\frac{(-1)^i}{i}=H_{\lfloor p/4\rfloor}-H_{(p-1)/2},\label{b-10}
\end{align}
and the following result (see \cite[(45)]{lehmer-am-1938}):
\begin{align}
H_{(p-1)/2}\equiv -2q_p(2)+pq_p(2)^2\pmod{p^2}.\label{b-11}
\end{align}

Finally, substituting \eqref{b-9} and \eqref{b-11} into \eqref{b-10}, we complete the proof of
\eqref{b-5}.
\qed

\section{Proof of Theorem \ref{t-1}}
By \eqref{a-cz}, we have
\begin{align}
\sum_{k=0}^{p-1}\frac{3k+1}{(-32)^k}{\rm Domb}(k)
&=\sum_{k=0}^{p-1}\frac{3k+1}{(-32)^k}\sum_{i=0}^k(-1)^i{k+2i\choose 3i}{2i\choose i}^2{3i\choose i}16^{k-i}\notag\\[10pt]
&=\sum_{i=0}^{p-1}\frac{1}{(-16)^i}{2i\choose i}^2{3i\choose i}\sum_{k=i}^{p-1}\frac{3k+1}{(-2)^k}{k+2i\choose 3i}.\label{c-1}
\end{align}
It can be easily proved by induction on $n$ that
\begin{align}
\sum_{k=i}^{n-1}\frac{3k+1}{(-2)^k}{k+2i\choose 3i}=(n-i){n+2i\choose 3i}(-2)^{1-n}.\label{c-2}
\end{align}
It follows from \eqref{c-1} and \eqref{c-2} that
\begin{align}
\sum_{k=0}^{p-1}\frac{3k+1}{(-32)^k}{\rm Domb}(k)
=\sum_{i=0}^{p-1}\frac{2^{1-p}(p-i)}{(-16)^i}{2i\choose i}^2{3i\choose i}
{p+2i\choose 3i}.\label{c-3}
\end{align}
Now we split the sum on the right-hand side of \eqref{c-3} into two pieces:
\begin{align*}
S_1=\sum_{i=0}^{(p-1)/2}(\cdot)\quad\text{and}\quad S_2=\sum_{i=(p+1)/2}^{p-1}(\cdot).
\end{align*}

For $0\le j\le (p-1)/2$, we have
\begin{align*}
(-1)^i(p-i){3i\choose i}{p+2i\choose 3i}
&=\frac{p(-1)^i(p+2i)\cdots(p+1)(p-1)\cdots(p-i)}{i!(2i)!}\\[10pt]
&=\frac{p(-1)^i(p+2i)\cdots(p+i+1)(p^2-1)\cdots(p^2-i^2)}{i!(2i)!}\\[10pt]
&\equiv \frac{pi!(p+2i)\cdots(p+i+1)}{(2i)!}\left(1-p^2H_i^{(2)}\right)\\[10pt]
&\equiv \frac{pi!(p+2i)\cdots(p+i+1)}{(2i)!}-p^3H_i^{(2)}\pmod{p^4}.
\end{align*}
Furthermore, we have
\begin{align*}
&\frac{pi!(p+2i)\cdots(p+i+1)}{(2i)!}\\[10pt]
&\equiv p\left(1+p\left(H_{2i}-H_i\right)+\frac{p^2}{2}\left((H_{2i}-H_i)^2-H_{2i}^{(2)}+H_i^{(2)}\right)\right)
\pmod{p^4}.
\end{align*}
It follows that
\begin{align*}
&(-1)^i(p-i){3i\choose i}{p+2i\choose 3i}\\[10pt]
&\equiv p+p^2\left(H_{2i}-H_i\right)+\frac{p^3}{2}\left((H_{2i}-H_i)^2-H_{2i}^{(2)}-H_i^{(2)}\right)
\pmod{p^4},
\end{align*}
and so
\begin{align}
S_1&\equiv 2^{1-p}p\sum_{i=0}^{(p-1)/2}\frac{1}{16^i}{2i\choose i}^2\notag\\[10pt]
&\times\left(1+p\left(H_{2i}-H_i\right)+\frac{p^2}{2}\left((H_{2i}-H_i)^2-H_{2i}^{(2)}-H_i^{(2)}\right)\right)
\pmod{p^4}.\label{c-4}
\end{align}

Note that for $0\le i\le \frac{p-1}{2}$,
\begin{align}
&(-1)^i{(p-1)/2\choose i}{(p-1)/2+i\choose i}\notag\\[10pt]
&=\frac{\left(\left(\frac{1}{2}\right)^2-\left(\frac{p}{2}\right)^2\right)
\left(\left(\frac{3}{2}\right)^2-\left(\frac{p}{2}\right)^2\right)\cdots
\left(\left(\frac{2i-1}{2}\right)^2-\left(\frac{p}{2}\right)^2\right)}{i!^2}\notag\\[10pt]
&\equiv \frac{1}{16^i}{2i\choose i}^2\pmod{p^2}.\label{c-5}
\end{align}
Letting $n=\frac{p-1}{2}$ in \eqref{b-1} and \eqref{b-2} and using \eqref{c-5}, we obtain
\begin{align}
&\sum_{i=0}^{(p-1)/2}\frac{1}{16^i}{2i\choose i}^2\left(H_{2i}-H_i\right)
\equiv (-1)^{\frac{p+1}{2}}\sum_{i=1}^{(p-1)/2}\frac{(-1)^i}{i}\pmod{p^2},\label{c-6}\\[10pt]
&\sum_{i=0}^{(p-1)/2}\frac{1}{16^i}{2i\choose i}^2
\left((H_{2i}-H_i)^2-H_{2i}^{(2)}-H_i^{(2)}\right)\notag\\[10pt]
&\equiv 2(-1)^{\frac{p-1}{2}}\left(\sum_{i=1}^{(p-1)/2}\frac{(-1)^i}{i^2}+\sum_{i=1}^{(p-1)/2}\frac{(-1)^i}{i}H_i\right)\pmod{p^2}.
\label{c-7}
\end{align}
Substituting \eqref{b-3}--\eqref{b-5} into the right-hand sides of \eqref{c-6} and \eqref{c-7} gives
\begin{align}
\sum_{i=0}^{(p-1)/2}\frac{1}{16^i}{2i\choose i}^2\left(H_{2i}-H_i\right)
\equiv (-1)^{\frac{p+1}{2}}\left(-q_p(2)+\frac{1}{2}pq_p(2)^2\right)+pE_{p-3}\pmod{p^2},
\label{c-8}
\end{align}
and
\begin{align}
&\sum_{i=0}^{(p-1)/2}\frac{1}{16^i}{2i\choose i}^2
\left((H_{2i}-H_i)^2-H_{2i}^{(2)}-H_i^{(2)}\right)\notag\\[10pt]
&\equiv (-1)^{\frac{p-1}{2}}q_p(2)^2+6E_{p-3}\pmod{p}.\label{c-9}
\end{align}
Moreover, by \cite[(1.7)]{sunzw-scm-2011} we have
\begin{align}
\sum_{i=0}^{(p-1)/2}\frac{1}{16^i}{2i\choose i}^2\equiv (-1)^{\frac{p-1}{2}}+p^2E_{p-3}\pmod{p^3}.
\label{c-10}
\end{align}
Substituting \eqref{c-8}--\eqref{c-10} into \eqref{c-4} and using the Fermat's little theorem, we arrive at
\begin{align}
S_1&\equiv (-1)^{\frac{p-1}{2}}p+5p^3E_{p-3}\pmod{p^4}.\label{c-11}
\end{align}

Next, we evaluate $S_2$ modulo $p^4$. For $(p+1)/2\le i \le p-1$, we have
${2i\choose i}^2\equiv 0\pmod{p^2}$, and
\begin{align*}
(-1)^i2^{1-p}(p-i){3i\choose i}{p+2i\choose 3i}&=\frac{(-1)^i2^{1-p}p(p+2i)\cdots(p+1)(p-1)\cdots(p-i)}{i!(2i)!}\\[10pt]
&\equiv \frac{p(p+1)\cdots(p+2i)}{(2i)!}\pmod{p^2}\\[10pt]
&=\frac{p(p+1)(p+2)\cdots 2p \cdots(p+2i)}{1\cdot2\cdots p\cdots 2i}\\[10pt]
&\equiv 2p\pmod{p^2},
\end{align*}
where we have utilized the Fermat's little theorem in the second step. Thus,
\begin{align*}
S_2\equiv 2p\sum_{i=(p+1)/2}^{p-1}\frac{1}{16^i}{2i\choose i}^2\pmod{p^4}.
\end{align*}
Recall the following supercongruence \cite[(1.9)]{sunzw-scm-2011}:
\begin{align*}
\sum_{i=(p+1)/2}^{p-1}\frac{1}{16^i}{2i\choose i}^2\equiv -2p^2E_{p-3}\pmod{p^3}.
\end{align*}
It follows that
\begin{align}
S_2\equiv -4p^3E_{p-3}\pmod{p^4}.\label{c-12}
\end{align}
Then the proof of \eqref{a-thm1} follows from \eqref{c-3}, \eqref{c-11} and \eqref{c-12}.

\section{Proof of Theorem \ref{t-2}}
By \eqref{a-sunzh}, we have
\begin{align}
\sum_{k=0}^{p-1}\frac{3k+2}{(-2)^k}{\rm Domb}(k)
&=\sum_{k=0}^{p-1}\frac{3k+2}{(-2)^k}\sum_{i=0}^{\lfloor k/2\rfloor}{2i\choose i}^2{3i\choose i}
{k+i\choose 3i}4^{k-2i}\notag\\[10pt]
&=\sum_{i=0}^{(p-1)/2}\frac{1}{16^i}{2i\choose i}^2{3i\choose i}\sum_{k=2i}^{p-1}
(-2)^k(3k+2){k+i\choose 3i}.\label{d-1}
\end{align}
Recall the following identity \cite[(2.4)]{guo-itsf-2013}:
\begin{align}
\sum_{k=2i}^{n-1}(-2)^k(3k+2){k+i\choose 3i}=(-1)^{n-1}(n-2i){n+i\choose 3i}2^n,\label{d-2}
\end{align}
which can be easily proved by induction on $n$.
It follows from \eqref{d-1} and \eqref{d-2} that
\begin{align}
\sum_{k=0}^{p-1}\frac{3k+2}{(-2)^k}{\rm Domb}(k)=
\sum_{i=0}^{(p-1)/2}\frac{2^p(p-2i)}{16^i}{2i\choose i}^2{3i\choose i}{p+i\choose 3i}.\label{d-3}
\end{align}

For $0\le i\le (p-1)/2$, we have
\begin{align*}
(p-2i){3i\choose i}{p+i\choose 3i}
&=\frac{p(p+i)\cdots(p+1)(p-1)\cdots(p-2i)}{i!(2i)!}\\[10pt]
&=\frac{p(p^2-1)(p^2-2^2)\cdots(p^2-i^2)(p-i-1)\cdots(p-2i)}{i!(2i)!}\\[10pt]
&\equiv \frac{p(-1)^ii!(p-i-1)\cdots(p-2i)}{(2i)!}\left(1-p^2H_i^{(2)}\right)\\[10pt]
&\equiv \frac{p(-1)^ii!(p-i-1)\cdots(p-2i)}{(2i)!}-p^3H_i^{(2)}\pmod{p^4}.
\end{align*}
Furthermore, we have
\begin{align*}
&\frac{(-1)^ii!(p-i-1)\cdots(p-2i)}{(2i)!}\\
&\equiv 1-p(H_{2i}-H_i)+\frac{p^2}{2}\left((H_{2i}-H_i)^2-H_{2i}^{(2)}+H_i^{(2)}\right)\pmod{p^3}.
\end{align*}
Thus,
\begin{align}
&(p-2i){3i\choose i}{p+i\choose 3i}\notag\\[10pt]
&\equiv p-p^2(H_{2i}-H_i)+\frac{p^3}{2}\left((H_{2i}-H_i)^2-H_{2i}^{(2)}-H_i^{(2)}\right)\pmod{p^4}.
\label{d-4}
\end{align}
Combining \eqref{d-3} and \eqref{d-4} gives
\begin{align}
\sum_{k=0}^{p-1}\frac{3k+2}{(-2)^k}{\rm Domb}(k)
&\equiv 2^p p\sum_{i=0}^{(p-1)/2}\frac{1}{16^i}{2i\choose i}^2\notag\\[10pt]
&\hskip-20mm\times\left(1-p(H_{2i}-H_i)+\frac{p^2}{2}\left((H_{2i}-H_i)^2-H_{2i}^{(2)}-H_i^{(2)}\right)\right)
\pmod{p^4}.\label{d-5}
\end{align}

Finally, substituting \eqref{c-8}--\eqref{c-10} into \eqref{d-5} and using the Fermat's little theorem, we obtain
\begin{align*}
\sum_{k=0}^{p-1}\frac{3k+2}{(-2)^k}{\rm Domb}(k)
&\equiv 2p(-1)^{\frac{p-1}{2}}\left(\left(2^{p-1}-1\right)^3+1\right)+6p^3E_{p-3}\\[10pt]
&\equiv 2p(-1)^{\frac{p-1}{2}}+6p^3E_{p-3}\pmod{p^4},
\end{align*}
as desired.

\section{Proof of Theorem \ref{t-3}}
By \eqref{a-ctyz}, we have
\begin{align*}
\sum_{k=0}^{n-1}(2k+1){\rm Domb}(k)8^{n-1-k}
&=\sum_{k=0}^{n-1}(2k+1)8^{n-1-i}
\sum_{i=0}^k(-1)^i{k\choose i}{k+i\choose i}f_i\\[10pt]
&=\sum_{i=0}^{n-1}(-1)^i8^{n-1-i}f_i\sum_{k=i}^{n-1}(2k+1)
{k\choose i}{k+i\choose i}.
\end{align*}
Note that
\begin{align*}
\sum_{k=i}^{n-1}(2k+1){k\choose i}{k+i\choose i}
=\frac{n(n-i)}{i+1}{2i\choose i}{n+i\choose 2i},
\end{align*}
which can be easily proved by induction on $n$. Thus,
\begin{align}
\frac{1}{n}\sum_{k=0}^{n-1}(2k+1){\rm Domb}(k)8^{n-1-k}
=\sum_{i=0}^{n-1}\frac{(-1)^i8^{n-1-i}(n-i)}{i+1}{2i\choose i}{n+i\choose 2i}f_i.
\label{e-1}
\end{align}
Since the Catalan numbers $C_i={2i\choose i}/(i+1)$ on the right-hand side of \eqref{e-1} are always integral, we conclude that the left-hand side of \eqref{e-1} is always a positive integer.

In a similar way, by using \eqref{a-ctyz} and the following identity:
\begin{align*}
\sum_{k=i}^{n-1}(-1)^k(2k+1){k\choose i}{k+i\choose i}
=(-1)^{n-1}n{n-1\choose i}{n+i\choose i},
\end{align*}
we obtain
\begin{align}
\frac{1}{n}\sum_{k=0}^{n-1}(2k+1){\rm Domb}(k)(-8)^{n-1-k}
=\sum_{i=0}^{n-1}(-1)^i8^{n-1-i}{n-1\choose i}{n+i\choose i}f_i.\label{e-2}
\end{align}
It is easy to see that the left-hand side of \eqref{e-2} is always an integer.

Next, we show that the left-hand side of \eqref{e-2} is positive.
From \cite[Proposition 2.8]{wz-scm-2014}, we conclude that the sequence $\{{\rm Domb}(k+1)/{\rm Domb}(k)\}_{k\ge 0}$
is strictly increasing. For $k\ge 2$, we have
\begin{align*}
\frac{{\rm Domb}(k+1)}{{\rm Domb}(k)}\ge \frac{{\rm Domb}(3)}{{\rm Domb}(2)}=\frac{64}{7}>8,
\end{align*}
and so the sequence $\{{\rm Domb}(k)/8^k\}_{k\ge 2}$ is strictly increasing. Let
\begin{align*}
a_k=\frac{(2k+1){\rm Domb}(k)}{8^k}.
\end{align*}
We immediately conclude that the sequence $\{a_k\}_{k\ge 0}$ is strictly increasing (the cases $k=0,1$ can be easily verified by hand).
Thus,
\begin{align*}
a_{n-1}-a_{n-2}+a_{n-3}-\cdots+(-1)^{n-1}a_0>0,
\end{align*}
and so
\begin{align*}
&\sum_{k=0}^{n-1}(2k+1){\rm Domb}(k)(-8)^{n-1-k}\\[5pt]
&=8^{n-1}\left(a_{n-1}-a_{n-2}+a_{n-3}-\cdots+(-1)^{n-1}a_0\right)>0.
\end{align*}
This proves the positivity for the left-hand side of \eqref{e-2}.

\vskip 5mm \noindent{\bf Acknowledgments.}
This work was supported by the National Natural Science Foundation of China (grant 11801417).

\end{document}